\documentclass[11pt,leqeq]{article}
\usepackage{amssymb,amsthm}
\usepackage{amsmath}
\usepackage{amscd}
\usepackage{xy}
\xyoption{all}
\usepackage[dvips]{graphicx}
\newtheorem{thm}{Theorem}
\newtheorem{lema}[thm]{Lemma}

\newtheorem{defi}[thm]{Definition}

\DeclareFontFamily{OMS}{rsfs}{\skewchar\font'60}
\DeclareFontShape{OMS}{rsfs}{m}{n}{<-5>rsfs5 <5-7>rsfs7 <7->rsfs10
}{} \DeclareSymbolFont{rsfs}{OMS}{rsfs}{m}{n}
\DeclareSymbolFontAlphabet{\scr}{rsfs}

\newcommand{\des}{\displaystyle}

\newcommand{\im}{{\rm{im}}}

\newcommand{\ma}{{\rm Map}}
\newcommand{\Co}{{\rm Cob}}

\newcommand{\Mfu}{{\rm MFunc}}

\newcommand{\Ho}{{\rm H}}

\newcommand{\Ca}{{\rm C}}

\newcommand{\ri}{{\rm r}}
\newcommand{\ve}{{\rm vect}}

\newcommand{\ob}{{\rm Ob}}
\newcommand{\HL}{{\rm HLQFT}}
\newcommand{\om}{{\rm Oman}}
\newcommand{\ca}{{\rm Cat}}
\newcommand{\op}{{\rm op}}

\renewcommand{\des}{\displaystyle}

\def\R{{\Bbb R}}

\def\R{{\mathbb R}}

\def\C{{\cal{C}}}

\def\d{{\cal{D}}}
\def\f{{\cal{F}}}
\def\g{{\cal{G}}}

\begin{document}
\title{Homology and manifolds with corners}
\author{Edmundo Castillo and
Rafael D\'\i az}
\maketitle

\begin{abstract}
We define a model for the homology of manifolds and use it to
describe the intersection product on the homology of compact
oriented manifolds and to define homological quantum field
theories which generalizes the notions of string topology
introduced by Chas and Sullivan and homotopy quantum field
theories introduced by Turaev.
\end{abstract}

\section{Introduction}

In this paper we propose a chain model for the homology groups of
smooth manifolds and show that this model is suited for the study
of the intersection product on the homology  of compact oriented
manifolds. Homology groups of a manifold $M$ are usually defined
trough the simplicial chain model $C^{s}(M)$  \footnote{All vector
spaces in this paper are defined over the complex numbers.} where
$$C^{s}(M)=\bigoplus_{i=0}^{\infty}C_{i}^{s}(M),$$
and $C_{i}^{s}(M)$ denotes the free vector space generated by
smooth maps from the $i$-dimensional simplex into $M$, i.e.,
$$C_{i}^{s}(M)=
\Big\langle c\ \ | \ \ c\colon \triangle_{i} \rightarrow M
\mbox{ is a smooth map }  \Big\rangle.$$

Above $\triangle_{i}= \{(x_0, \dots ,x_i) \in  (\mathbb{R}_{\geq
0})^{i+1}\ \ | \ \ x_{0} + ... + x_{i} = 1
\}.$  For $0 \leq k \leq i$ we have maps $e_{k}\colon \triangle_{i-1} \rightarrow \triangle_{i}$
given by $e_{k}(x_0,...,x_{i-1})=
(x_0,...,x_{k-1},0,x_k,...,x_{i-1})$. The differential
$\partial\colon C_{i}^{s}(M) \rightarrow C_{i-1}^{s}(M)$ is given by
$\partial(c)=
\sum_{k=0}^{i}(-1)^{k}c \circ e_k$ for all $\ c:
\triangle_{i} \rightarrow M.$  The singular homology $H(M)$ of the manifold $M$ are by definition
the homology groups of the complex $C^{s}(M)$. \\

We are interested in the intersection product on homology. Despite
its simplicity and good properties the simplicial chains model
have major shortcomings when it comes to define  the intersection
product at the chain level. The main difficulties can be
summarized in the following facts
\begin{itemize}
\item{Given chains $c: \bigtriangleup_i \rightarrow M$ and $d: \bigtriangleup_j \rightarrow M$
one obtains a map $c \times d:\bigtriangleup_i \times
\bigtriangleup_j
\rightarrow M \times M$ which unfortunately does not define a chain on $M \times M$ since
$\bigtriangleup_i \times \bigtriangleup_j$ is a not a simplex. The
way around is to triangulate $\bigtriangleup_i \times
\bigtriangleup_j$ which can be done in a canonical but not unique way.
The choice of triangulations introduces a combinatorics somewhat
foreign to the geometric setting on which homology  should lie.}

\item{The natural domain for the intersection product  $c\cap d$ of \textit{transversal} chains
$c$ and $d$ is $(c\times d)^{-1}(\Delta_M)$, where $\Delta_M =
\{(m,m) \in M \times M \ \ | \ \ m \in M
\}$. Even if we assume as given  canonical triangulations for
$\bigtriangleup_i \times \bigtriangleup_j$ this will be of little
help in order to define triangulations on $(c\times
d)^{-1}(\Delta_M)$. Such triangulations exist but are neither
unique nor canonical.}
\end{itemize}

These facts explain why it is so hard to define the intersection
product at the chain level using the simplicial chain model. The
main problem lies in the severe restrictions on the domain of
simplicial chains, indeed only standard simplices are allowed. In
order to overcome this difficulties we propose an alternative
model for the homology of manifolds which allows a wide a variety
of domains for its chains. For each manifold we define a complex
$C(M)=\bigoplus_{i=0}^{\infty}C_{i}(M)$. The space of degree $i$
chains $C_{i}(M)$ is constructed as follows
\begin{itemize}
\item{Let $\overline{C}_{i}(M)$ be the vector space freely generated by pairs
$(K,c)$ such that $K$ is a compact oriented manifold with corners
and $c:K \rightarrow M$ is a smooth map.}
\item{$C_{i}(M)$ is the quotient of $\overline{C}_{i}(M)$ by the
relations\\

1)$(K^{\op},c) = -(K,c)$ where $K^{\op}$ is the manifold $K$
provided with the opposite orientation.\\

2) $(K_1 \sqcup K_2, c_1 \sqcup c_2)= (K_1,c_1) +(K_2,c_2)$.}
\end{itemize}

The  differential $\partial \colon \Ca_{i}(M) \to
\Ca_{{i}-\rm {1}}(M)$ is given by
$$\partial(K_,c)=\sum_{L \in \pi_{0}(\partial_{1}
K)}(\overline{L},c|_{ \overline{L}}).$$

The sum above ranges over the connected components $L$ of the
first boundary strata $\partial_{1} K$ of $K$ provided with the
induced boundary orientation. $c|_{\overline{L}}$ denotes the
restriction of $c$ to the closure of $L.$ Complexes $C(M)$ enjoy
the following properties that make them worthy of study

\begin{itemize}
\item{The homology groups $\Ho(M)$ of the complex $C(M)$ are isomorphic to
the singular homolo\-gy groups of $M$. The isomorphism map is
induced by the inclusion $i\colon C(M)  \rightarrow C(M),$ which is well
defined since every simplex is a manifold with corners.}
\item{Given chains $c\colon K \rightarrow M$ and $d:L \rightarrow M$ we have a
chain $c\times d\colon K \times L \rightarrow M \times M$, since the
Cartesian product of oriented manifolds with corners is an
oriented manifold with corner. Moreover one can show that
$$\partial_{1}( K
\times L )= (\partial_{1} K)\times L \ \
+ \ \ (-1)^{dim(K)}K \times (\partial_{1} L).$$}
\item{If  $c:K \rightarrow M$ and $d:L \rightarrow M$ are
\textit{transversal} chains then  $(c \times d)^{-1}(\Delta_M)$ is a submanifold with
corners of $K \times L$. Letting $\pi: K\times L$ be the
projection onto $K$, the intersection product $c\cap d$ of
\textit{transversal} chains $c$ and $d$ is defined to be
$((c \times d)^{-1}(\Delta_M)\  , \  c \circ \pi)$. }

\end{itemize}

The purpose of this paper is to systematically use the manifold
with corners chain model  to study the intersection product on the
homology of compact connected oriented manifolds and its various
generalizations. The paper is organized as follows

\begin{itemize}
\item{In Section \ref{sec2} we formalize our construction of the
complex $C(M)$ for each manifold $M$. We show that
homology of $C(M)$ agrees with the singular homology of  $M$. }

\item{In Section \ref{sec3} we extend the notion of \textit{transversality}
for maps whose domain is a manifold with corners and whose target
is a smooth manifold. We use this construction to give a
characterization of the intersection product on the homology
groups of compact oriented manifolds.}

\item{In Section \ref{sec4} we use our model for homology to show that
the operad $\Ho(D_d)$ of little discs in dimension $d$ acts on the
homology groups $\Ho(M^{S^{d}})$ of the space of free $d$-spheres on
a compact connected oriented manifold $M$, thus obtaining a new
proof of the following result due to Voronov
\cite{V1}, see also the book by  Cohen, Hess, and Voronov
\cite{CoVo}.

{\noindent}{\bf Theorem
\ref{te2}.} For $d\geq 1$, the graded vector space $\Ho(M^{S^{d}})$  is

\begin{enumerate}
\item  A differential graded associative algebra if $d=1.$
\item A differential graded twisted Poisson algebra with the commutative
associative product of the degree $0$ and with the Lie bracket of
degree $(1-d)$ in the case of odd $d\geq 3.$
\item A differential graded twisted Gerstenhaber algebra with the commutative
associative product of the degree $0$ and with the Lie bracket of
degree $(1-d)$ in the case of even $d\geq 2.$
\end{enumerate}}

\item{In Section \ref{sec5} we construct the category $\Ho(M^{S(Y)})$ of dynamical $Y$ branes
living on a manifold $M$ and show that this category admits
analogues of the notion of transposition and trace for matrices.}

\item{In Section \ref{sec7} we apply our model for homology to study homological quantum field
theories \HL \ first introduced in our works \cite{Cas2} and
\cite{Cas1}. The notion of homological quantum field theory \HL
\ which lies in the crossroad of two lines of thought. On one hand we have the notion
of topological quantum field theory as axiomatized by Atiyah in
\cite{MA} and further generalized  to the notion of homotopical
quantum field theory by Turaev in \cite{TU1} and \cite{TU2}. On
the other hand we have the already mention loop product on
$\Ho(M^{S^{1}})$ defined by Chas and Sullivan in \cite{SCh} and
further studied by Cohen and Jones \cite{CJ}, and Cohen, Jones and
Jun \cite{CJJ}. Combining both constructions we arrived to the
notion of
\HL \ . We give an example of homological
quantum field theory in arbitrary dimension.}
\end{itemize}

\section{Homology using manifolds with corners}\label{sec2}

Let us star this Section with a brief introduction to manifolds
with corners. For $0 \leq k \leq n$, we denote by $\ H^{n}_{k}$
the subspace of $\mathbb{R}^{n}$ given by
\[H^{n}_{k}=[0,\infty)^{k}\times \mathbb{R}^{n-k}=\{(x_{1},\cdots
,x_{n})\in \mathbb{R}^{n} \ \ | \ \ x_{1} \geq 0, \cdots ,x_{k}
\geq 0
\}.\]A subset $V \subset H^{n}_{k}$ is open if there exists an open
subset $W \subset \mathbb{R}^{n}$  such that $V=W\cap H^{n}_{k}.$
Let $\partial_{0}(H^{n}_{k})=(0,\infty)^{k}\times
\mathbb{R}^{n-k},$ and for a set $V\subseteq
H^{n}_{k},$  let $\partial_{0}(V)=V\cap \partial_{0}(H^{n}_{k}).$
A map $f\colon V \to \R$ is said to be smooth if there exits open
set $W\subseteq\R^{n}$ and smooth map $F\colon W \subseteq \R^{n}
\to \R$ such that $V=W\cap H^{n}_{k}$ and $F\mid_{V}=f.$ Given open sets
$V_{i}\subset H^{n}_{k_{i}}, \ 0\leq k_{i}\leq n,
\ i=1,2,$ we say that a map $f \colon V_{1}\to V_{2}$ is a
diffeomorphism if it is a homeomorphism with inverse $g
\colon V_{2} \to V_{1}$ and each coordinate component of $f$ or $g$
is a smooth map.\\

Let $M$ be an Hausdorff topological space. $M$ is a $n$-manifold
with corners if it is locally homeomorphic to $H^{n}_{k}$ for some
$0\leq k\leq n,$ that is, there exists an open cover
$\mathfrak{U}=\{U_{i}\}_{i \in \Lambda}$ of $M$ such that for each
$i\in \Lambda,$ there is a map $\varphi_{i}\colon U_{i} \to
H^{n}_{k_{i}},\ 0\leq k_{i}\leq n,$ which maps $U_{i}$
homeomorphically onto an open subset of $H^{n}_{k_{i}}.$ We call
$(\varphi_{i},U_{i})$ a chart with domain $U_{i}.$ The set of
charts $\Phi=\{(\varphi_{i},U_{i})\}_{i \in \Lambda}$ is an atlas.
Two charts $(\varphi_{i},U_{i}),(\varphi_{j},U_{j})$ are said to
have a smooth overlap if the coordinate changes
\[ \varphi_{j}\circ \varphi_{i}^{-1}\colon \varphi_{i}(U_{i}\cap U_{j})
  \to \varphi_{j}(U_{i}\cap U_{j} ) \]
are smooth diffeomorphisms. An atlas $\Phi$ on $M$ is called
smooth if every pair of charts in it have smooth overlaps. There
is a unique maximal smooth atlas which contains $\Phi.$ A maximal
atlas $\Phi$ on $M$ defines a  structure of a smooth manifold with
corners on $M.$ The pair $(M,\Phi )$ is called a $n$-manifold with
corners.\\

All manifold with corner are naturally stratified spaces with
smooth strata. The smooth strata of $H^{n}_{k}$ are given for
$0\leq l\leq k$ by
\[\partial_{l}H^{n}_{k}=\{x\in H^{n}_{k}\mid x_{i}=0 \
\mbox{for exactly}\ l \ \mbox{of the first}\
  k \ \mbox{indices}\}.\]
Notice that
\[\txt \footnotesize {$ \partial_{l}H^{n}_{k}=
\des\bigsqcup_{\substack{  I\subseteq \{1,\cdots ,k\}\\
\mid
 I\mid=l}} H^{n}_{I}$}\]
where $H^{n}_{I}=\{(x_{1},\cdots ,x_{n})\mid x_{i}=0  \mbox{ if
and only if } i\in I \}.$ For a  manifold with corners $M$ we set
\[\partial_{l} M=\{m\in M\mid \ \mbox{there exists local
coordinates mapping } m \
 \mbox{ to} \ \partial_{l}H^{n}_{k} \}.\]
One can checks that $ M= \des
\bigsqcup_{0\leq l \leq n} \partial _{l} M$ and that each
$\partial _{l} M$ is a submanifold of $M.$ A manifold with corners
$M$ is a  manifold if $M =\partial_{0}M.$ A manifold with corners
is a manifold with boundaries if $\partial_{2}M=
\emptyset.$  Given a manifold $M$ we define the graded vector
space
\[\Ca(M)=\bigoplus_{i=0}^{\infty} \Ca_{i}(M),\]
where  $\Ca_{\it i}(M)$ denotes the
 vector space
\[\frac{\Big\langle (K,c)\colon \begin{array}{c}
K\ \mbox{is an  oriented}\ i \mbox{-manifold with }   \\
\mbox{corners and}  \ c\colon K \to M \
\mbox{is a smoth map} \\
\end{array}
\Big\rangle}{\langle (K^{\op},c)-(K,c), (K\sqcup L,c \sqcup d) - (K,c) - (L,d) \rangle}.\]

\begin{figure}[ht]
\begin{center}
\includegraphics[height=2.3cm]{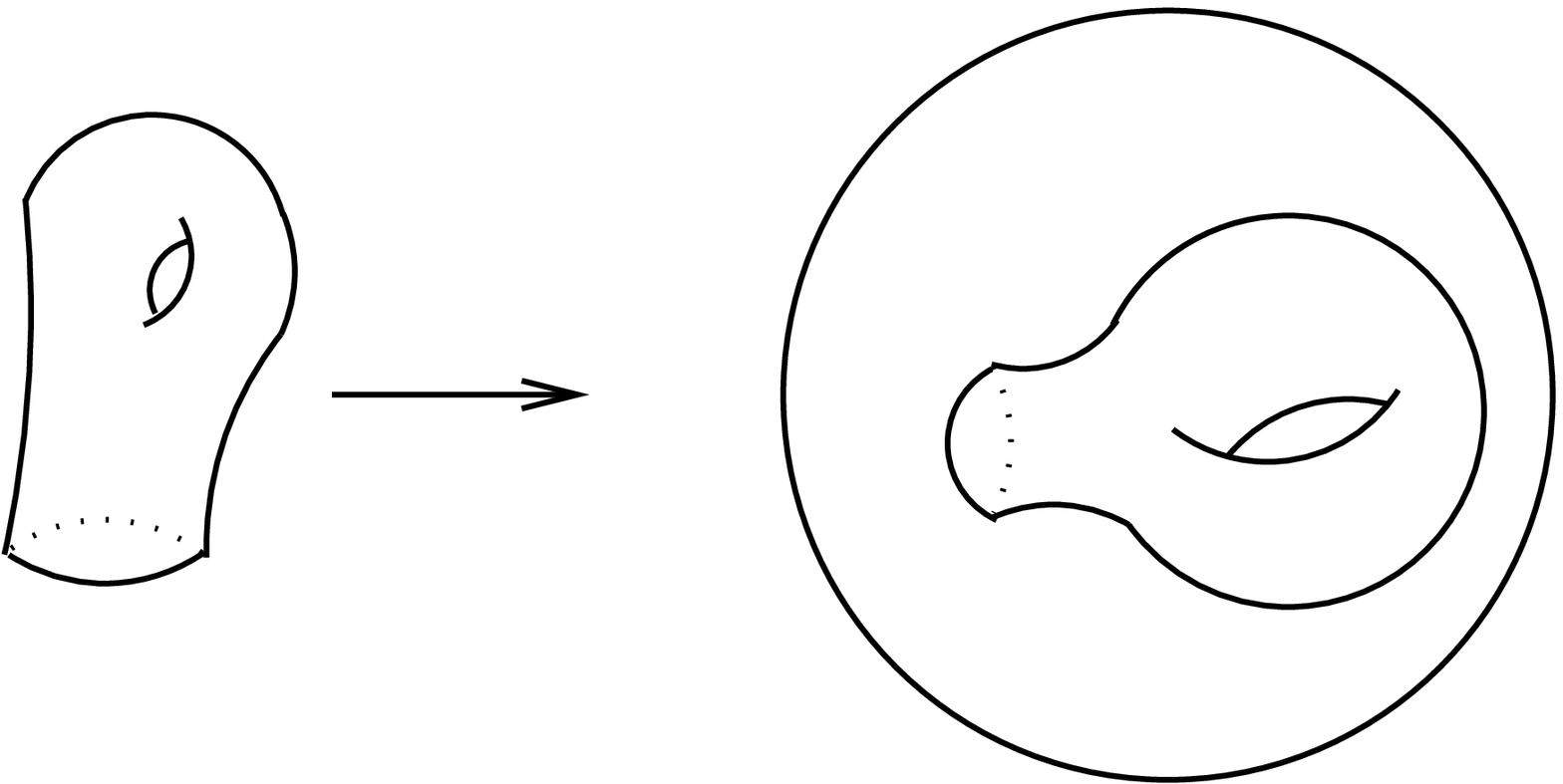}
\caption{\ Chain with domain a manifold with corners.\label{fig:Ge1}}
\end{center}
\end{figure}

\begin{figure}[ht]
\begin{center}
\includegraphics[height=2.3cm]{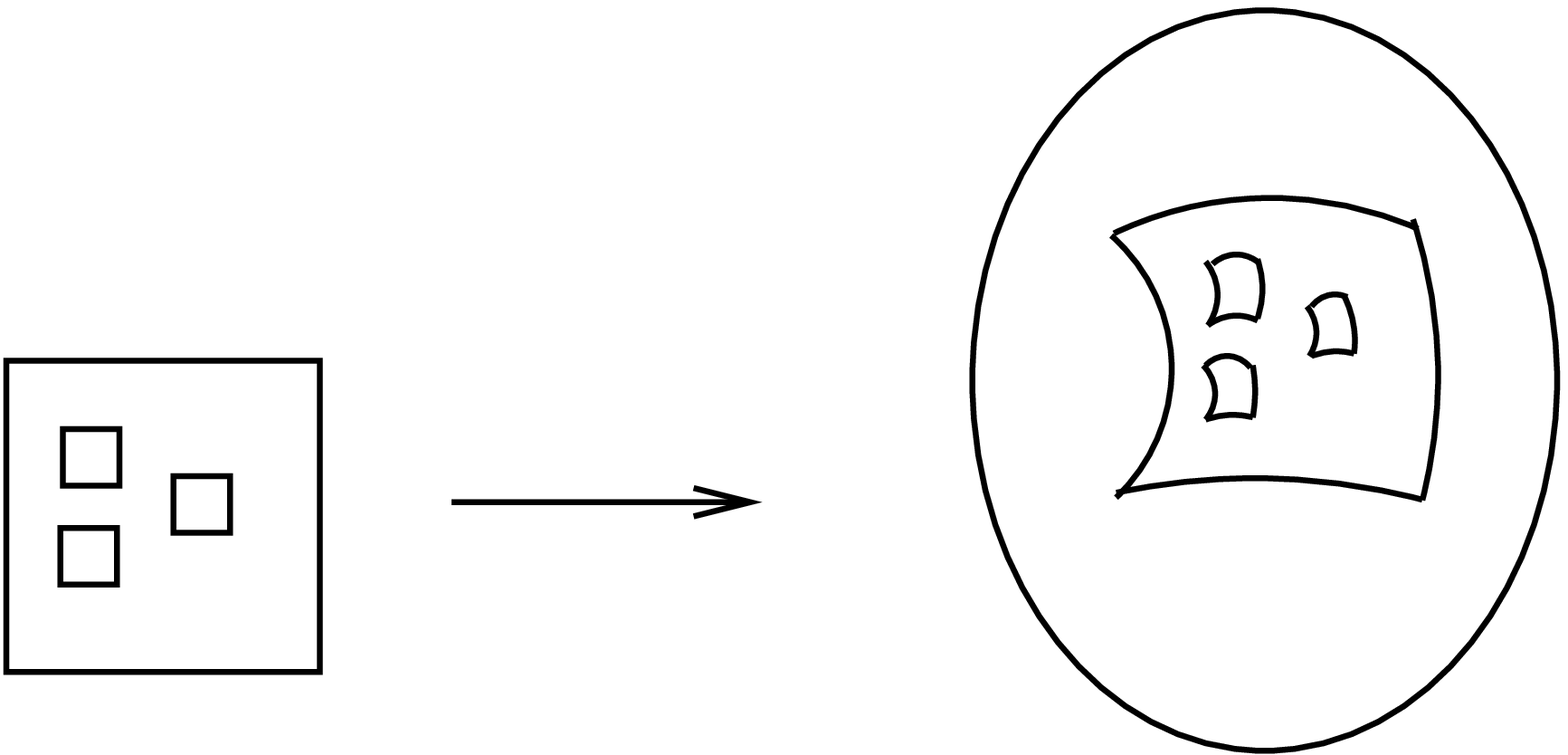}
\caption{\ Chain with domain a manifold with corners.
\label{fig:cua1}}
\end{center}
\end{figure}

We recall that $M^{\op}$ denotes the manifold $M$ provided with the
opposite orientation. Figures \ref{fig:Ge1} and \ref{fig:cua1}
below show a couple of examples of chains whose domain are
manifold with corners.

\begin{defi}
 {\em Let the map $\partial \colon \Ca_{i}(M) \to \Ca_{{i}-\rm {1}}(M)$
be given by $\partial(K,c)=\sum_{L}(\overline{L},
c|_{\overline{L}}),$ where the sum runs over the connected
components of $\partial_{1} K$ provided with the induced boundary
orientation and $c|_{\overline{L}}$ denotes the restriction of $c$
to the closure of $L.$}
\end{defi}

With this notation we have the following

\begin{thm}
 {\em $(\Ca(M),\partial)$ is a differential complex. Moreover}
 $\Ho(\Ca(M),\partial)=\Ho(M)$.
\end{thm}

\begin{proof}
$\partial^{2}=0$ since

$$\partial^{2}(K,c)= \sum_{L \in \pi_0(\partial_{2} K) }[(\overline{L},
c|_{\overline{L}}) + (\overline{L}^{\op}, c|_{\overline{L}})] .$$

There is an obvious inclusion $i\colon C^{s}(M)\to
\Ca(M).$ The map $i$ is quasi-isomorphism since any manifold with
corners can be triangulated \footnote{We thank J. P. Brasselet, M.
Goresky, and  R. Melrose for helpful comments on the triangulation
of manifold with corners. See references \cite{MG}, \cite{JEF}
 and \cite{ve} for more information on triangulability.} and thus any chain in $\Ca(M)$ is homologous to a chain
in $C^s(M).$
\end{proof}

\section{Transversal intersection of manifolds with corners}\label{sec3}

In this Section we study  transversal intersection for manifold
with corners. Let us first state a useful fact.
\begin{lema}\label{Le9}
 {\em Let $M$ and $N$ be manifolds with corners then $M\times N$ is
also a manifold with corners. Moreover if $x\in
\partial_{l} M$ and $y\in
\partial_{k} N$ then  $(x,y)\in
\partial_{l+k} (M\times N).$}
\end{lema}

\begin{defi}
 {\em Let $N_{0} \ (N_{1})$ be a $n_{0}\ (n_{1})$-manifold with
corners and $M$ be a smooth oriented manifold. Let $f_{0}\colon
N_{0} \to M$ and $f_{1}\colon N_{1} \to M$ be smooth maps. We say
that $f_{0}$ and $f_{1}$ are {\it transversal,}  $f_{0}\pitchfork
f_{1},$ if for $0\leq k\leq n_{0}, \ 0\leq s\leq n_{1}, \
f_{0}\mid_{\partial_{k}(N_{0})}$ and
$f_{1}\mid_{\partial_{s}(N_{1})}$ are transversal, i.e., for
$x_{0}\in\partial_{k}(N_{0})$ and $x_{1}\in \partial_{s} (N_{1})$
such that $f_{0}(x_{0})=f_{1}(x_{1})=m$ we have that
$d(f_{0})(T_{x_{0}}\partial_{k}(N_{0}))
+d(f_{1})(T_{x_{1}}\partial_{s}(N_{1}))=T_{m}M.$}
\end{defi}

The next couple of results are generalization of the corresponding
results for manifold with boundaries proved by Guillemin and
Pollack in
\cite{GP}.
\begin{lema}\label{2}
 {\em Let $M$ be a smooth manifold and $f=(f_{1},\cdots f_{s})\colon
M \to \R^{s}$  a smooth map with regular value $(0,\cdots ,0).$
Then $\{x\in M \mid f_{i}(x)\geq 0, \ i=1,\cdots ,s \}$ is a
manifold with corners.}
\end{lema}

\begin{proof}
$\{x \in M \mid f(x)> 0\}$ is a open subspace of $M$ thus it is a
smooth manifold. Let $x\in M$ be such that $f(x)=0.$ Since $0$ is
a regular value then $f$ is locally equivalent to the submersion
$\pi^{n}_{s}\colon \R^{n}\to \R^{s}$ given by $\pi^{n}_{s}
(x_{1},\cdots , x_{n})= (x_{1},\cdots , x_{s}).$ The desired
result holds since locally $\{x\in M\mid f_{i}(x)\geq 0, \
i=1,\cdots ,s\}= \{x\in
\R^{n} \mid x_{i}\geq 0, \ i=1,\cdots ,s\}.$
\end{proof}

\begin{thm}\label{1}
 {\em Let $M$ be a  $m$-manifold and $N$ be a $n$-manifold with
corners. Let $f\colon N \to M$ be a smooth map and $i\colon
P\hookrightarrow M$ be an embedded $p$-submanifold of $M$. If
$f\pitchfork i$ then $f^{-1}(P)$ is a submanifold with corners of
$N.$}
\end{thm}

\begin{proof}
Since $f\mid_{\partial_{0}(N)}$ is transversal to $i$ we have that
$f^{-1}(P)\cap \partial_{0} (N)$ is a manifold of codimension
$m-p.$ Consider $x\in f^{-1}(P)\cap\partial_{k} (N)$ and  let $l$
be a submersion from a neighborhood of  $f(x)$ in $M$ to
$\mathbb{R}^{m-p}$ such that in this neighborhood $P=l^{-1}(0).$
 Choose a local parametrization $h\colon U \to N$ around
$x$ where $U$ is an open subset of $H^{n}_{t},$ $0\leq k \leq t
\leq n,$ and set $g=l\circ f\circ h.$ Since $h$ is a
diffeomorphism $f^{-1}(P)$ is a manifold with corners in a
neighborhood of $x$ if and only if $(f\circ h)^{-1}(P)=g^{-1}(0)$
is a manifold with corners near $y=h^{-1}(x)\in \partial_{k} (U).$
The tranversality condition
$d(f)(T_{x}\partial_{k}(N))+T_{f(x)}(P)=T_{f(x)}(M)$ implies that
$x$ is a regular point of $l\circ f,$ i.e., that $g$ is regular at
$y.$ Since $g$ is smooth then it can be  extended to a smooth map
$G$ on a neighborhood  of $y$ in $\mathbb{R}^{n}.$ As $dG=dg, \ G$
is also regular at $y$ and thus $G^{-1}(0)$ is a smooth
submanifold of $\mathbb{R}^{n}$ in a neighborhood of $y.$\\

  Since $g^{-1}(0)=G^{-1}(0)\cap H^{n}_{t}$ in a neighborhood of $y$, we
must show that $G^{-1}(0)\cap H^{n}_{t}$ is a manifold with
corners in a neighborhood of $y.$ The map $\pi=(\pi_{1},\cdots
,\pi_{t})\colon G^{-1}(0)\subseteq \R^{n} \to
\R^{t},$ is such that  $G^{-1}(0)\cap H^{n}_{t}= \{s\in G^{-1}(0) \mid
\pi_{i} (s)\geq 0, \ i=1,\cdots ,t\}.$\\

We claim that $(0,\cdots ,0)$ is a regular value for $\pi.$
Otherwise,  there is a point $z\in G^{-1}(0)$ such that
$\pi(z)=(0,\cdots ,0)$ and $\dim (d_{z}\pi)< t.$ Thus, there exist
at least a linear relation $a_{1}d_{z}\pi_{1}+\cdots
+a_{t}d_{z}\pi_{t}=0.$ Making a further change of variables we may
assume that the linear relation is $d_{z}\pi_{1}=0.$\\

The fact that $d_{z}\pi_{1}$ is zero on $T_{z}(G^{-1}(0))$ means
that the first coordinate of every vector in $T_{z}(G^{-1}(0))$ is
zero, i.e., that $T_{z}(G^{-1}(0))\subset
\{0\}\times\R^{n-1}.$ The kernel of $dg=dG\colon \R^{n} \to \R$ is
$T_{z}(G^{-1}(0))$ and $d_{z}(\partial_{\{1\}} g)$ is the
restriction of $d_{z}g\colon \R^{n} \to \R$ to $\R^{n-1}.$ Since
the kernel of $dg$ is contained in $\R^{n-1}$ then the linear maps
$dg\colon \R^{n}\to \R^{1}$ and $d(\partial_{\{1\}} g)\colon
\R^{n-1}\to \R$ must have the same kernel. By transversality
both maps are surjective, so the kernel of $dg$ has dimension
$n-1$ whereas the kernel of $d(\partial_{\{1\}} g)$ has dimension
$n-2,$ which is a contradiction. Lemma  \ref{2} finishes the proof
of Theorem \ref{1}.
\end{proof}

\begin{lema}
 {\em Let $K$ and $L$ be manifolds with
corners and $M$ be an oriented manifold. Assume that $c\colon K
\to M$ and $d\colon L\to M$ are transversal maps, then $$K
\times _{M}L= \{ (a,b)\in K\times L\mid x(a)=y(b)
\}$$ is an oriented manifold with corners. Moreover
$T(K\times_{M}L)=TK\times_{TM}TL.$}
\end{lema}

\begin{proof}
 $K\times L$ is a manifold with corners By Lemma \ref{Le9}. Since
$c$ and $d$ are transversal the map $c\times d\colon K\times L\to
M\times M$ is  transversal to $\Delta
\colon M
\to M\times M.$   Theorem \ref{1} then guarantees that $(c\times
d)^{-1}(\Delta(M))=K\times_{M} L$ is a manifold with corners.
\end{proof}

\begin{lema}
 {\em Let $K$ and $L$ be manifolds with
corners and $M$ be a oriented manifold. Assume also that $c\colon
K \to M$ and $d\colon L\to M$ are transversal maps. Then

\begin{itemize}\label{seguridad}
\item{$\partial_{1}( K
\times L )= (\partial_{1} K)\times L \ \
+ \ \ (-1)^{dim(K)}K \times (\partial_{1} L)$.}

\item{$\partial_{1}( K \times _{M}L) = (\partial_{1} K)\times_{M} L \ \
+ \ \ (-1)^{dim(K)}K \times_{M} (\partial_{1} L)$.}

\end{itemize}}
\end{lema}

Lemma \ref{seguridad} guarantees that the intersection product on
homology $\Ho(M)$ to be defined below trough the intersection of
transversal chains is indeed well defined. The following result
gives us information on the algebraic structure on the homology of
a compact connected oriented manifold $M$.

\begin{thm}
 {\em  $\Ho(M)[dim(M)]$ is a Frobenius algebra.}
\end{thm}

The Frobenius algebra structure on $\Ho(M)[dim(M)]$ is given

\begin{itemize}
\item{An intersection product $\cap\colon \Ho(M)[dim(M)] \otimes \Ho(M)[dim(M)] \rightarrow \Ho(M)[dim(M)]$,
given by $[K,c]\cap [L,d] = [K \times_{M} L, c \circ \pi]$ for
transversal chains $c: K \rightarrow M$ and $d: L
\rightarrow M$.}
\item{There is a trace map $tr\colon \Ho(M)[dim(M)] \rightarrow \mathbb{C}$ defined to be zero in non-zero
degrees and the identity in degree $-d$.}
\item{The bilinear form $< \ \, \ \ >\colon \Ho(M)[dim(M)] \otimes \Ho(M)[dim(M)] \rightarrow \mathbb{C}$ defined by
$<a,b>=tr(a \cap b)$, for $a,b \in \Ho(M)[dim(M)],$ is
non-degenerated by Poincar\'e duality. }
\item{ The unit in $\Ho(M)[dim(M)]$ is defined by the identity map  $I\colon M \rightarrow M$.}

\end{itemize}

\section{Algebraic structures on $\Ho(M^{S^{d}})[dim(M)]$ }\label{sec4}

Let $D^{d}=\{ x\in \mathbb{R}^{d} \mid \|x \|\leq 1 \}.$ A little
disc in $D^{d}$ is an affine transformation $T_{a,r}\colon
D^{d}\to D^{d}$ given by
 $T_{a,r}(x)=rx+a,$  where $0 < r < 1$ and
$a \in D^{d}.$ For $n\geq 0,$ let us introduce the topological
spaces
\[D_{d}(n)= \Big\{(T_{x_{1},r_{1}}, \ldots ,T_{x_{n},r_{n}}) \ \Big |
\begin{array}{c}
x_{i}, x_{j}\in D^{d},\ 0\leq r_{i} < 1 \ \mbox{such that if }\
i\neq j\
 \\
 \overline{\im(T_{x_{i},r_{i}})}\cap \overline{\im(T_{x_{j},r_{j}})}=
\emptyset \ \mbox{ for all} \ 1\leq
i,j \leq n  \\
\end{array} \Big \} .\]

{\noindent}Notice that the disc with center $a$ and radius $r$ is
obtained as the image of the transformation $T_{a,r}$ applied the
standard disc $D^{d}$. The little discs operad   was introduced by
J. M. Boardman and R. M. Vogt  \cite{BV}, and Peter May
\cite{May1}.

\begin{defi}
 {\em The {\it little discs operad} in dimension $d,$
$D_{d}=\{D_{d}(n)\}_{n\geq 0}$ is the topological operad with the
following structure}
\begin{enumerate}

\item  {\em The composition $\gamma_{k}:D_{d}(k)\times D_{d}(n_{1})\times
\cdots \times D_{d}(n_{k}) \to D_{d}(j),$ with $j=\sum n_{s},$ are
defined as follows: given $(T_{x_{1},r_{1}}, \cdots
,T_{x_{k},r_{k}}) \in D_{d}(k)$ and $(T_{x_{i1},r_{i1}}, \cdots
,T_{x_{in_{i}},r_{in_{i}}})\in D_{d}(n_{i})$ for $1\leq i\leq k$,
by $$\gamma_{k}\big((T_{x_{1},r_{1}}, \cdots
,T_{x_{k},r_{k}});(T_{x_{11},r_{11}}, \cdots
,T_{x_{1n_{1}},r_{1n_{1}}}),\cdots ,(T_{x_{k1},r_{k1}}, \cdots ,
T_{x_{kn_{k}},r_{kn_{k}}})\big)= $$ $$(T_{x_{1},r_{1}}\circ
T_{x_{11},r_{11}}, \cdots ,T_{x_{1},r_{1}}\circ
T_{x_{1n_{1}},r_{1n_{1}}},\cdots , T_{x_{k},r_{k}}\circ
T_{x_{k1},r_{k1}}, \cdots ,T_{x_{k},r_{k}}\circ
T_{x_{kn_{k}},r_{kn_{k}}}).$$}

\item  {\em The action of the symmetric group $S_{n}$ on
$D_{d}(n)$ is given by  $$\sigma(T_{x_{1},r_{1}},\cdots ,
T_{x_{n},r_{n}})=(T_{x_{\sigma^{-1}(1)},r_{\sigma^{-1}(1)}},\cdots
, T_{x_{\sigma^{-1}(n)},r_{\sigma^{-1}(n)}}),$$ for all $\sigma
\in S_{n}$ and $(T_{x_{1},r_{1}},\cdots , T_{x_{n},r_{n}})\in
D_{d}(n).$}
\end{enumerate}
\end{defi}

Since $D_{d}$ is a topological operad its homology  $\Ho(D_{d})$ is
an operad in the category of graded vector spaces. So the
following definition makes sense

\begin{defi}\label{d-algebras}
 {\em A $d$-algebra is an algebra over $\Ho(D_{d})$ in the category of graded vector
 spaces.}
\end{defi}

Let us fix a compact oriented manifold $M$.

\begin{defi}\label{def13}
 {\em We denote by $M^{S^{d}}$ the set of all smooth maps $\alpha
\colon D^{d}\to M$ such that $\alpha$ is  constant in an open
neighborhood of $\partial( D^{d}).$ $M^{S^{d}}$ is given the
compact-open topology.}
\end{defi}

Let $
  \Ca(M^{S^{d}})=\des\bigoplus_{i=0}^{\infty}
\Ca_{i}(M^{S^{d}})$ be the vector space generated by
chains $c\colon K \to M^{S^{d}}$ such that the associated map
$\widehat{c}\colon K\times D^{d} \to M$ given by
$\widehat{c}(k,p)=(c(k))(p)$ for $(k,p)\in K\times D^{d}$ is
smooth. Let $e\colon M^{S^{d}}\to M$ be the map given by
$e(\alpha)=\alpha(\partial (D^{d})).$ We denote by $e_{*}$ the
induced maps $e_{*}\colon
\Ca(M^{S^{d}}) \to \Ca(M)$ and $e_{*}\colon \Ho(M^{S^{d}}) \to \Ho(M)$.

\begin{thm}\label{te3}
 {\em $\Ho(M^{S^{d}})[dim(M)]$ is a $d$-algebra.}
\end{thm}
\begin{proof}We construct  maps $\theta_{n}\colon \Ho(D_{d}(n)) \otimes
\Ho(M^{S^{d}})^{n} \to \Ho(M^{S^{d}})$. Assume that $[K,c]
\in
\Ho(D_{d}(n))$ and $ [K_{i},c_i] \in \Ho(M^{S^{d}})$ for $1\leq i
\leq n.$ Let $[L,d]= \theta_{n}([K,c];[K_{1},c_1],\ldots ,[K_{n},c_n])$ be defined as follows
\begin{itemize}
\item{$L= K \times K_1 \times_M ... \times_M K_n$.}

\item{In order to define $d:L \rightarrow  M^{S^{d}}$ let $(k,k_1,...,k_n) \in L$
and assume that $c(k)$ is such that
$c(k)=(T_{p_{1}(k),r_{1}(k)},\cdots , T_{p_{n}(k),r_{n}(k)}),$
then the map $d(k;k_{1},\ldots ,k_{n})\colon S^{d} \longrightarrow
M$ is given for $y \in S^{d}$ by }
\end{itemize}

\[ d(k;k_{1},\ldots ,k_{n})(y) =
 \left\{ \begin{array}{ll}
e(c_{1}(k_{1}))  \ \ \mbox{if}\ y
\notin \bigcup \overline{\im(T_{p_{i}(k),r_{i}(k)})}\\
  \ &  \\
c_{i}(k_{i})\Big(\frac{y-p_{i}(k)}{r_{i}(k)}\Big) \ \ \mbox{if}\ y
\in \overline{\im(T_{p_{i}(k),r_{i}(k)})}
\end{array} \right. \]
\end{proof}

Next result follows from Theorem \ref{te3} above and Theorem 3 of
\cite{KO1}.

\begin{thm}\label{te2}
 {\em  For $d\geq 1,$ $\Ho(M^{S^{d}})[dim(M)]$ is}
\begin{enumerate}
\item  {\em A differential graded associative algebra if $d=1.$}
\item  {\em A differential graded twisted Poisson algebra with the
commutative associative pro\-duct of the degree $0$ and with the
Lie bracket of degree $(1-d)$ in the case of odd $d\geq 3.$ \item
A differential graded twisted Gerstenhaber algebra with the
commutative associative product of degree $0$ and with the Lie
bracket of degree $(1-d)$ in the case of even $d\geq 2.$}
\end{enumerate}
\end{thm}

\section{Branes category}\label{sec5}

Let $M$ be a compact connected oriented manifold. Let $N_{0}$ and
$N_{1}$ be connected oriented embedded submanifolds of $M.$ Let
$Y$ be a compact manifold. \\

\begin{defi}
{\em Let $M^{S(Y)}(N_{0},N_{1})$ be the set of smooth maps $f\colon
Y\times [-1,1] \to M$ such that}
\begin{itemize}
\item{$f(y,-1)\in N_{0},$ $f(y,1)\in N_{1}$ for $y \in Y$.}
\item{\em $f$ is constant on open neighborhoods of $Y\times\{-1\}$ and
$Y\times\{1\}$, respectively.}
\end{itemize}
\end{defi}

We give $M^{S(Y)}(N_{0},N_{1})$ the compact-open topology. Notice
that $M^{S(Y)}(N_{0},N_{1})$ is a subspace of $\ma(S(Y),M)$ where
$S(Y)= Y\times[-1,1]/\sim$ and $\sim$ is the equivalence relation
given by $y_{1}\times
\{-1\}\sim y_{2}\times\{-1\},$ and $y_{1}\times
\{1\}\sim y_{2}\times\{1\}$ for  $y_{1},y_{2}\in Y.$\\

Let  $\Ca(M^{S(Y)}(N_{0},N_{1}))=\txt
\footnotesize {$\des\bigoplus_{i=0}^{\infty}$}
\Ca_{i}(M^{S(Y)}(N_{0},N_{1}))$ be the vector space generated
by chains $c\colon K\to M^{S(Y)}(N_{0},N_{1})$  such that the maps
$$e_{-1}(x)\colon K\times Y\to N_{0} \mbox{ and } e_{1}(x)\colon
K\times Y\to N_{1}$$ given by $e_{i}(c)(x,y)=c(x)(y,i)$ for
$i=-1,1$ are smooth. Let $$e_{-1}\colon M^{S(Y)}(N_{0},\newline
N_{1})\to N_{0} \mbox{ and }  e_{1}\colon M^{S(Y)}(N_{0},N_{1})\to
N_{1}$$ be the maps given by $e_{-1}(f)=f(y,-1)\in N_{0}$ and
$e_{1}(f)=f(y,1)\in N_{1},$ respectively. We denote by $e_{i*}$
the induced maps $$e_{i*}\colon
\Ca(M^{S(Y)}(N_{0},N_{1})) \to \Ca(N_{0}) \mbox{ and } e_{i*}\colon
\Ho(M^{S(Y)}(N_{0},N_{1})) \to \Ho(N_{0})$$ for $i=-1,1.$\\

Let us introduce the category $\Ho(M^{S(Y)})$ given by

\begin{itemize}
\item{Objects in $\Ho(M^{S(Y)})$ are oriented connected compact embedded submanifolds of
$M$.}

\item{Morphisms in $\Ho(M^{S(Y)})$  are
 given by $\Ho(M^{S(Y)})(N_0, N_1)= \Ho(M^{S(Y)}(N_0, N_1)).$}

\item{The identity morphism $[N,I_N] \in \Ho(M^{S(Y)})(N, N)$ is given by the map
$$\overline{I}_N\colon N \times S(Y) \longrightarrow M,$$
given by $\overline{I}_N(n,s)=n$ for $n \in N$ and $s \in S(Y)$.}
\item{ Given $[K,c] \in \Ho(M^{S(Y)})(N_0, N_1)$ and $[L,d] \in \Ho(M^{S(Y)})(N_1, N_2)$ the composition
morphism $[L,d] \circ [K,c] \in \Ho(M^{S(Y)})(N_0, N_1)$ is given
for transversal cycles $c$ and $d$ by the map
$$\widehat{[L,d] \circ [K,c]}\colon L \times_M K \times Y \times [-1,1] \longrightarrow M$$
defined by

\[\widehat{[L,d] \circ [K,c]}(l,k,y,t)
 = \left\{ \begin{array}{ll}
\widehat{c}(k,y,t)  \ \ \mbox{if} \ t \in [0,\frac{1}{2}] \\
 \ & \ \\
\widehat{d}(l,y,t) \ \ \mbox{if}\
t \in [\frac{1}{2},1]

          \ & \ \\
\end{array} \right. \]}
\end{itemize}

Summing up we obtain

\begin{thm}
{\em $\Ho(M^{S(Y)})$ is a category.}
\end{thm}

 An interesting feature of the category
 $\Ho(M^{S(Y)})$ is that it comes equipped with analogues of the operation of transposition
 and trace on matrices.

 \begin{thm}
{\em There is a contravariant functor $\ri\colon \Ho(M^{S(Y)})\to
\Ho(M^{S(Y)})$ identical on objects, given for $N_{0}$ and $N_{1}$
by
\[\ri \colon \Ho(M^{S(Y)})(N_{0},N_{1})\to
\Ho(M^{S(Y)})(N_{1},N_{0})\] taking $[K,c]$ into $[K,\ri(c)]$
defined as follows: for $k \in K$, $y \in Y$ and  $-1\leq t
\leq 1,$ we have that
$$\ri(c)(k)](y,t)=c(k)(y,-t).$$}
\end{thm}

Let $M^{Y\times S^{1}}$ be the space of smooth maps $f\colon
Y\times S^{1} \to M$ provided with the compact-open topology. We
denote by $\Ho(M^{Y\times S^{1}})=\des \bigoplus_{i=0}^{\infty}
\Ho_{i}(M^{Y\times S^{1}})$ the vector space generated by chains
$c\colon K\to M^{Y\times S^{1}}$ such that the induced map
$\widehat{c}\colon K\times Y\times S^{1}\to M$ is smooth.

\begin{thm}
{\em For each object $N$ in $\Ho(M^{Y\times S^{1}})$  there is a map
$$tr\colon \Ho(M^{S(Y)})(N,N)
\longrightarrow \Ho(M^{Y\times S^{1}}),$$ where for $[K,c]$ such that
$e_{-1 *}(c)$ is transversal to $e_{1 *}(c)$ we define
$$tr([K,c]= [((e_{-1 *}(c) \times e_{1 *}(c)) \circ \Delta)^{-1}(\Delta_M) , c \circ \pi ],$$
where $\Delta\colon M \longrightarrow M \times M$ is the diagonal map.}
\end{thm}

\section{Homological quantum field theory}\label{sec7}

\qquad The theory  of cobordism was introduced by Rene Thom
\cite{RT}. Based on the notion of cobordisms Michael F. Atiyah
\cite{MA} and G. B. Segal \cite{Se} (see also \cite{Baez},
\cite{Kim}, and \cite{Toen}) introduced the axioms for topological
quantum field theories TQFT and conformal field theories CFT,
respectively. V. Turaev  \cite{TU1} introduced the axioms for
homotopical quantum field theories. Essentially the axioms of
Atiyah may be summarized as follows: 1) Considerer the monoidal
category $\Co_{n}$ of $n$-dimensional cobordism. 2) Define the
category of TQFT as the category $\Mfu(\Co_{n},\ve)$ of monoidal
functors $\f\colon \Co_{n}\to \ve.$ Segal's definition of CFT may
also be  recast as a category of monoidal functors.\\

\begin{figure}[ht]
\begin{center}
\includegraphics[height=2.3cm]{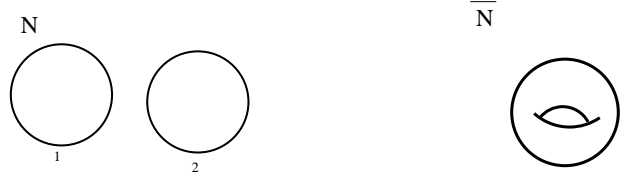}
\caption{ \ Manifold $N$ and its completion
$\overline{N}.$\label{fieti}}
\end{center}
\end{figure}

Let us fix a few conventions. $\pi_{0}(M)$ denotes the set of
connected components of $M.$  We define the completion of a
manifold $N$ to be $\overline{N}=\des\prod_{c\in
\pi_{0}(N)}c.$ Figure \ref{fieti} gives an example of a manifold
and its completion.
We denote by $D(M)$ the set of embedded connected oriented
submanifolds of $M.$  By convention the empty set is assumed to be
a $n$-dimensional manifold for all $n\in
\mathbb{N}.$
\\

Following the pattern discussed above we define the category
$\HL_d(M)$ of homological quantum field theories as follows

\begin{enumerate}

\item{First we construct a category $\Co_{d}^{M},$ object of which are
called extended cobordisms.}

\item{ Second we define $\HL_d(M)$ to be the category of
monoidal functors $\Mfu(\Co_{d}^{M},\ve).$}

\end{enumerate}

Objects in $\Co_{d}^{M}$ are triples $(N,f,<)$ such that

\begin{itemize}
\item{$N$ is a compact oriented  manifold of dimension $d-1$.}
\item{$f\colon \pi_{0}(N) \to D(M)$  is a map.}
\item{$<$ is a linear ordering on
 $\pi_{0}(N).$ .}
\end{itemize}

For objects $(N_{0},f_{0},<_{0})$ and  $(N_{1}, \ f_{1},<_{1})$ in
$\Co_{d}^{M}$ we set
\[ \Co_{d}^{M}((N_{0},f_{0},<_{0}), \ (N_{1},f_{1},<_{1}))=
\overline{\Co_{d}^{M}}((N_{0},f_{0},<_{0}), \ (N_{1},f_{1},<_{1}))
\diagup \backsim\] where $\overline{\Co_{d}^{M}}((N_{0},f_{0},<_{0}), \
(N_{1},f_{1},<_{1}))$
is the set of triples $(P,\alpha,[c])$ such that\\
\begin{enumerate}
\item{$P$ is a compact oriented $d$-manifold with boundaries.}

\item{ $\alpha \colon N_{0}\bigsqcup N_{1}\times [0,1)\to
\im (\alpha)\subseteq P$ is a diffeomorphism.  $\alpha \mid_{
N_{0}\bigsqcup N_{1}}\to
\partial  P$ is such that $\alpha|_{N_{0}}$
reverses the orientation and $\alpha|_{N_{1}}$ preserves the
orientation.}

\item{ $[c] \in \Ho(\ma (P,M)_{f_{0}, f_{1}})[dim(\overline{N_{1}})],$ where $\ma(P,M)_{f_{0},
f_{1}}$ denotes the set of smooth maps $g\colon P\to M$ such that
$g$ is constant on a neighborhood of each connected component of
its boundary $\partial  P.$}
\end{enumerate}

If $\alpha\in
\ma(P,M)_{f_{0},f_{1}}$ we define $e_{i}(\alpha)\in \overline{N_i}$ by
$e_{i}(\alpha)(t)=\alpha(t)$
for $t\in \pi_{0}(N_i)$ and $i=0,1.$\\

 Triples $(P,\alpha,\xi )$ and $(P',\alpha', \xi' )$
in $\overline{\Co_{d}^{M}}((N_{0},f_{0},<_{0}), \
(N_{1},f_{1},<_{1}))$  are equivalent according to  $\sim$ if and
only if

\begin{itemize}
\item{ There exists an orientation
preserving diffeomorphism $\varphi\colon P_{1} \to P_{2}$ such
that:}
\item{ $\varphi \circ
\alpha=\alpha'$.}
\item{$\varphi_{\star}(\xi)=\xi'.$}
\end{itemize}

Assume we are given $$(P_{1},\alpha_{1},[c_{1}])\in
 \Co_{d}^{M}
((N_{0},f_{0},<_{0}), (N_{1}, f_{1}, <_{1}))$$ and
$$(P_{2},\alpha_{2},[c_{2}])\in
 \Co_{d}^{M}
((N_{1},f_{1},<_{1}), (N_{2}, f_{2}, <_{2})),$$

Assume $[c] \in \Ho(\ma (P,M)_{f_{0}, f_{1}})$ is given by a smooth
map $c
\colon K\times P\to M$ for a compact oriented manifold with
corners $K.$ \\

 The composition
morphism $$(P_{2},\alpha_{2},[c_2])\circ
(P_{1},\alpha_{1},[c_{1}]))\in
\Co_{d}^{M}((N_{0},f_{0},<_{0}),
(N_{2},f_{2},<_{2})))$$  is the triple $(P_{2}\circ
P_{1},\alpha_{2} \circ
\alpha_{1},[c_{2} \circ c_{1}])$ where

\begin{itemize}
\item{ $ P_{2}\circ P_{1}=P_{1}\des\bigsqcup_{N_{1}} P_{2}$.}

\item{$\alpha_{2} \circ
\alpha_{1}=\alpha_{2}\mid_{N_{2}}\bigsqcup
\alpha_{1}\mid_{N_{0}} $.}
          \\
\item{ The domain of $c_{2} \circ c_{1}$ is $L \times_{\overline{N}_1}K$. The map
$c_{2} \circ c_{1}:L \times_{\overline{N}_1}K \rightarrow \ma
(P_{1}\circ P_{2},M)_{f_{0}, f_{2}}$ has associated map
$$\widehat{c_{1} \circ c_{2}}:L \times_{\overline{N}_1}K \times (P_{1}\circ
P_{2}) \rightarrow M$$ given by $\widehat{c_{2} \circ
c_{1}}(l,k,p)= \widehat{c_1}(k,p)$ if $p
\in P_1$ and $\widehat{c_{2} \circ c_{1}}(k,l,p)=
\widehat{c_2}(l,p)$ if $p \in P_2.$ Figure $\ref{Bor3}$ represents a $n$-cobordism enriched over $M$
and Figure $\ref{Bor4}$ shows a composition of $n$-cobordism
enriched over $M.$
\begin{figure}[h]
\begin{minipage}[t]{0.5\linewidth}
\begin{center}
\includegraphics[height=3.5cm]{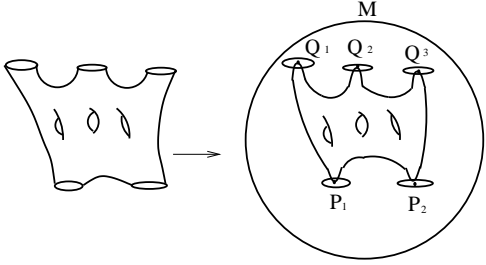} \caption{ \ A
$n$-cobordism enriched over $M.$}\label{Bor3}
\end{center}
\end{minipage}
\begin{minipage}[t]{0.5\linewidth}
\begin{center}
\includegraphics[height=4.2cm]{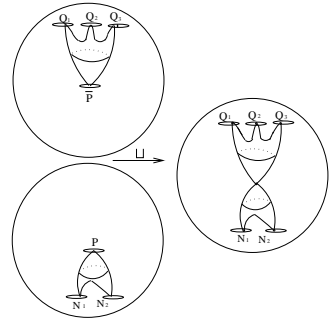} \caption{ \ Composition of
$n$-cobordisms enriched over $M.$}\label{Bor4}
\end{center}
\end{minipage}
\end{figure} }

\item{The identity morphism $(N \times [0,1] ,\alpha,[c])\in
 \Co_{d}^{M}( (N,f,<), (N, f, <))$ is determined by the map
 $\widehat{c}:N \times N \times [0,1] \longrightarrow M$ given by
 $\widehat{c}(n,m,t)=n$ for $n,m \in N$ and $t \in [0,1]$.   }

\end{itemize}

Thus we obtain

\begin{thm}
 {\em $(\Co_{n}^{M}, \ \sqcup, \ \emptyset)$ is a monoidal category
with disjoint union $\sqcup$ as product and empty set as unit. }
\end{thm}

\begin{defi}
 {\em $(\Co_{n,r}^{M}, \ \sqcup, \ \emptyset)$ is the full monoidal subcategory
of $\Co_{n}^{M}$ without the unit. }
\end{defi}

Given monoidal categories $\C$ and $\d$ we denote by $\Mfu(\C,\d)$
the category of monoidal functors from $\C$ to $\d.$ We are
finally ready to introduce our main definition.

\begin{defi}
{\em $\HL_{d}(M)=\Mfu(\Co_{n}^{M},\ve)$ and $\HL_{d,r}
(M)=\Mfu(\Co_{d,r}^{M},\ve).$
 $\HL_{d}(M)$ is the category of $d$ dimensional homological quantum field theories.
 $\HL_{d,r}(M)$  is the category of the $d$ dimensional restricted homological quantum field theories.}
\end{defi}

Explicitly the category $\HL_{d}(M)$ is given by

\begin{itemize}
\item{Objects in $\HL_d(M)$ are monoidal functors $\f \colon \Co_{n}^{M}
\to \ve.$}

\item{Morphisms in $\HL_d(M)(\f ,\g)$ are natural transformations  $T\colon \f \to
\g $ for objects  $\f, \g$ in $\HL_d(M).$}
\end{itemize}

 We now give a construction which yields our main example of homological field theory.
 Consider the map $\Ho \colon \Co_{d,r}^{M} \to
\mbox{vect}$ given on objects by
\[\xymatrix @C=.1in  @R=.02in
{\Ho \colon \ob(\Co_{d,r}^{M})
 \ar[rr] & &
\ob(\mbox{vect})\\
(N,f,<) \ar @{|->}[rr] & & \
 \Ho(N,f,<)= \Ho(\overline{N})
 }\]

The image under $\Ho$ of a morphism $(P,\alpha, [c])\in
\Co_{n,r}^{M}((N_{0},f_{0},<_{0}), (N_{1}, \ f_{1},<_{1}))$  is
the linear map
\[\xymatrix @C=.1in  @R=.02in {\Ho(P,\alpha,[c])
    \colon \Ho(N_{0},f_{0},<_{0}) \ar[rr] & &
\Ho(N_{1},f_{1},<_{1})\\
\ \ \ \ \ \ \ \ \ \ \ \ \ \ \ [h] \ar @{|->}[r] & &
[ \Ho(P, \alpha,\xi)(h)]}  \]

 defined as follows:\\

Let $[d] \in
\Ho(\ma(P,M)_{f_{0},f_{1}})$
 be given by the chain $d \colon L\to \ma(P,M)_{f_{0},
f_{1}}.$  Below we use the map $e_{0}(d)
\colon L \to \overline{N}_{0}.$   Assume that  $h$ is given by a map
$h: O \rightarrow \overline{N}_{0}$. We define the domain of
$\Ho(P, \alpha ,d)(h)$  by $L \times_{\overline{N}_{0}} O$ and let
the  map $$\Ho(P, \alpha ,d)(h):L
\times_{\overline{N}_{0}} O  \rightarrow  \overline{N}_1$$ be given by

$$\Ho(P, \alpha ,d)(h)(l,o)=e_{1}(d(l)),$$ for
$(l,o) \in L
\times_{\overline{N}_{0}} O.$

\begin{thm}
 {\em $\Ho$ defines a restricted homological quantum field theory.}
\end{thm}
\begin{proof}
Let $$(P_{1},\alpha_{1}, [c_{1}])\in
\Co_{n,r}^{M}((N_{0},f_{0},<_{0}), (N_{1}, \ f_{1},<_{1}))$$ and
$$(P_{2},\alpha_{2}, [c_{2}])\in \Co_{n,r}^{M}((N_{1},f_{1},<_{1}),
(N_{2}, \ f_{2},<_{2})),$$ where $c_1 : K_1 \rightarrow
\ma(P,M)_{f_{0},f_{1}}$ and $c_2 : K_2 \rightarrow
\ma(P,M)_{f_{1},f_{2}}.$  We want to check that
\[\Ho((P_{2},\alpha_{2}, [c_{2}])\circ (P_{1},\alpha_{1},
[c_{1}]))=\Ho(P_{2},\alpha_{2}, [c_{2}])\circ\Ho(P_{1},\alpha_{1},
[c_{1}]).\] The domain of $(P_{2},\alpha_{2}, c_{2})\circ
(P_{1},\alpha_{1}, c_{1})$ is
$K_{2}\times_{\overline{N}_{1}}K_{1},$ thus for $[h] \in
H(\overline{N}_0)$ given by a chain $h:L \rightarrow
\overline{N}_0$ the domain of
$$\Ho((P_{2},\alpha_{2}, c_{2})\circ
(P_{1},\alpha_{1}, c_{1}))(h)$$ is
$(K_{2}\times_{\overline{N}_{1}}K_{})\times_{\overline{N}_{0}}L.$
On the other hand the domain of $$\Ho(P_{2},\alpha_{2}, c_{2})(
\Ho(P_{1},\alpha_{1},
c_{1})(h))$$ is
$K_2\times_{\overline{N}_{1}}(K_{1}\times_{\overline{N}_{0}}L).$
Thus we see that the domains agree. It is easy to check that the
correponding maps also agree.
\end{proof}
\vspace{0.3cm}

We close this paper studying the functorial properties of the
correspondence $$M
\to \HL_d(M).$$ Let \om \ be the groupoid whose objects are
compact oriented smooth manifolds and whose morphisms are
orientation preserving diffeomorphisms. By
\ca \ we denote the category of small categories.

\begin{thm}
{\em The map $\HL_d(-)\colon\om^{\op} \longrightarrow \ca$ given
on objects by $M\to \HL_d(M)$ is functorial.}
\end{thm}
\begin{proof}
Let $M,N$ be objects in $\om$ and $\varphi\in \om^{\op}(M,N)$ an
orientation preserving diffeomorphism $\varphi\colon N
\to M$. $\varphi$ induces a functor
 $\varphi_{\star}\colon \Co_{d}^{N} \to \Co_{d}^{M}$ which is given on
 objects by
 $\varphi_{\star}(N,f,<)=(N,\varphi\circ f, \varphi_{*}(<))$ and similarly on morphisms.
Thus we get a functor $\varphi^{\star}\colon \Mfu(\Co_{d}^{M},\ve)
\to \Mfu(\Co_{d}^{N},\ve)$ given by $\varphi^{\star}(\f)=\f\circ
\varphi_{\star}$  for any monoidal functor
 $\f\colon  \Co_{d}^{N} \to \ve$.
\end{proof}

\vspace{0.3cm}
\subsection*{Acknowledgment} Thanks to Nicolas Andruskiewistch, Delia Flores de Chela,
Takashi Kimura, Bernardo Uribe, Sylvie Paycha and Arturo Reyes.

\bibliographystyle{amsplain}
\bibliography{homologywithcorners}

\providecommand{\bysame}{\leavevmode\hbox to3em{\hrulefill}\thinspace}
\providecommand{\MR}{\relax\ifhmode\unskip\space\fi MR }
\providecommand{\MRhref}[2]{%
  \href{http://www.ams.org/mathscinet-getitem?mr=#1}{#2}
}
\providecommand{\href}[2]{#2}
\begin{thebibliography}{10}

\bibitem{MA}
M.~F. Atiyah, \emph{Topological quantum field theory}, Publications
  math\'ematiques de {L' I.H.\'E.S}. 17 (1982), no.~4, 661--692.

\bibitem{Baez}
J.~C. Baez and J.~Dolan, \emph{Higher-dimensional {Algebra} and {Topological}
  {Quantum} {Field} {Theory}}, J. Math. Phys. 36 (1995), 6073--6105.

\bibitem{BV}
J.~M. Boardman and R.~M. Vogt, \emph{Homotopy invariant algebraic structures on
  topological spaces}, vol. 347, Springer Verlag, 1973.

\bibitem{Cas2}
E.~Castillo and R.~D\'{\i}az, \emph{Homological matrices}, math.KT/0510443, to
  appear in Contemporary Mathematics, 2005.

\bibitem{Cas1}
\bysame, \emph{Homological quantum field theory}, math.KT/0509532, 2005.

\bibitem{SCh}
M.~Chas and D.~Sullivan, \emph{String {Topology}}, Math.GT/9911159, too appear
  in Annals of Mathematics, 1999.

\bibitem{CoVo}
R.~L. Cohen, K.~Hess, and A.~Voronov, \emph{String topology and cyclic
  homology}, vol. Adv. Courses Math. CRM Barcelona, Birkh{\"a}user, Basel,
  2006.

\bibitem{CJ}
R.~L. Cohen and J.~D.~S. Jones, \emph{A homotopy theoretic realization of
  string topology}, Math. Ann. 324 (2002), no.~4.

\bibitem{CJJ}
R.~L. Cohen, J.~D.~S. Jones, and J.~Yan, \emph{The loop homology algebra of
  spheres and projective spaces. {C}ategorical decomposition techniques in
  algebraic topology}, Progr. Math., (2004), no.~215.

\bibitem{MG}
M.~Goresky, \emph{Triangulation of stratified objects}, Proc. Amer. Math. Soc.
  72 (1978), no.~1, 193--200.

\bibitem{GP}
V.~Guillemin and A.~Pollack, \emph{Differential {Topology}}, Prentice-Hall,
  1974.

\bibitem{JEF}
F.~E.~A. Johnson, \emph{On the triangulation of smooth fibre bundles}, Fund.
  Math. 118 (1983), no.~1, 39--58.

\bibitem{Kim}
T.~Kimura, \emph{Topological {Quantum} {Field} {Theory} and {Algebraic}
  {Structures}}, in Quantum field theory and noncommutative geometry workshop
  at Tohoku University in Sendai, 2003.

\bibitem{KO1}
M.~Kontsevich, \emph{Operads and motives in deformation quantization}, Letter
  in Mathematical Physics (1999), no.~48, 3572.

\bibitem{May1}
J.~P. May, \emph{The geometry of {Iterated} {Loop} {Spaces}}, vol. Lecture
  Notes in Mathematics 271, Springer Verlag, 1972.

\bibitem{Se}
G.~B. Segal, \emph{Othe definition of conformal field theory}, Topology,
  geometry and quantum field theory, London Math. Soc. Lecture Note Ser., 308
  (2004), 421--576.

\bibitem{RT}
R.~Thom, \emph{Sur les vari\'et\'es cobordantes.}, Colloque de topologie et
  géométrie différentielle, {Strasbourg} (1952), no.~7.

\bibitem{Toen}
B.~Toen, \emph{Notes on higher categorical structures in topological quantum
  field theory}, Lecture in MPI-Bonn, June 2000.
  http://www.picard.ups-tlse.fr/~toen/note.html, 2000.

\bibitem{TU1}
V.~Turaev, \emph{Homotopy field theory in dimension {2} and group-algebras},
  QA/9910010, 1999.

\bibitem{TU2}
\bysame, \emph{Homotopy field theory in dimension 3 and crossed
  group-categories}, math.GT/0005291, 2000.

\bibitem{ve}
A.~Verona, \emph{Stratified mappings---structure and triangulability}, vol.
  Lecture Notes in Mathematics, {\textbf1102}, Springer Verlag, 1984.

\bibitem{V1}
A.~Voronov, \emph{Notes on universal algebra. {G}raphs and {P}atterns in
  {M}athematics and {T}heoretical {P}hysics ({M}. {L}yubich and {L}.
  {T}akhtajan, eds.)}, Proc. Sympos. Pure Math. 73. AMS, Providence, RI (2005),
  81--103.

\end{thebibliography}

\noindent Edmundo Castillo. \\
ecastill@euler.ciens.ucv.ve \\
Universidad Central de Venezuela\\

\noindent Rafael D\'\i az.\\
ragadiaz@gmail.com \\

\end{document}